\RequirePackage[leqno]{amsmath} 
\documentclass[leqno]{amsart}
\usepackage{caption}
\usepackage[utf8]{inputenc}     
\usepackage{ifthen} 
\usepackage{multicol}
\usepackage{multirow}
\usepackage{floatrow}
\usepackage{longtable}
\usepackage[singlelinecheck=false 
]{caption}
\usepackage{color, graphics}
\usepackage{graphicx, multicol, latexsym, amsmath, amssymb}
\usepackage{makecell}
\setcellgapes{4pt}
\usepackage{blindtext}
\usepackage{subfigure}
\usepackage{amsfonts}
\usepackage{enumitem}
\setlist[itemize]{topsep=0pt,after=\vspace{1.5\baselineskip}} 
\usepackage[colorlinks=true]{hyperref}
\usepackage{empheq}
\usepackage[T1]{fontenc}
\usepackage{tabularx}
\usepackage{pgfplots}
\usepackage{tikz}
\tikzset{  
mynode/.style={fill,circle,inner sep=2pt,outer sep=0pt}
}
\usetikzlibrary{arrows.meta,arrows}
\usepackage[margin=0.7in]{geometry}
\usetikzlibrary{backgrounds,automata}

\setlist[itemize]{noitemsep, topsep=0pt}

\def\R{\mathbb R} \def\N{\mathbb N}

\def\R{\mathbb R} \def\N{\mathbb N} 
\def\TM{T_{\rm{max}}} 
\makeatletter
\@namedef{subjclassname@2020}{\textup{2020} Mathematics Subject Classification}
\makeatother

\def
\@cite
#1#2{[{{\bfseries #1}\if@tempswa , #2\fi}]}
\newtheorem{theorem}{Theorem}[section]

\newtheorem{lemma}[theorem]{Lemma}

\newtheorem{remark}{Remark}

\title[Uniformly bounded solutions to an attraction-repulsion
model] 
      {
Improvements and generalizations of results concerning attraction-repulsion chemotaxis models
}
\author[S. Frassu, T. Li and G. Viglialoro]{}
\subjclass[2020]{Primary: 35A01, 35K55, 35Q92. Secondary:  92C17.}
\keywords{Chemotaxis, Global existence, Boundedness, Nonlinear production, Consumption, Logistic source. \\
\textit{$^*$Corresponding author}: silvia.frassu@unica.it}

\begin{document}

\maketitle

\centerline{\scshape{Silvia Frassu$^{\sharp,*}$, Tongxing Li$^{\natural}$ \and Giuseppe Viglialoro$^{\sharp}$}}
\medskip
{
\medskip
\centerline{$^\sharp$Dipartimento di Matematica e Informatica}
\centerline{Universit\`{a} di Cagliari}
\centerline{Via Ospedale 72, 09124. Cagliari (Italy)}
\medskip
}
{
\medskip
\centerline{$^{\natural}$School of Control Science and Engineering} 
\centerline{Shandong University}
\centerline{Jinan, Shandong, 250061 (P. R. China)}
\medskip
}
\bigskip
\begin{abstract}
We enter the details of two recent articles concerning as many chemotaxis models, one nonlinear and the other linear, and both with produced chemoattractant and saturated  chemorepellent.  More precisely, we are referring respectively to the papers ``Boundedness in a nonlinear attraction-repulsion Keller--Segel system with production and consumption'', by S. Frassu, C. van der Mee and G. Viglialoro [\textit{J. Math.\ Anal.\ Appl.} {\bf 504}(2):125428, 2021] and ``Boundedness in a chemotaxis system with consumed chemoattractant and produced chemorepellent'', by S. Frassu and G. Viglialoro [\textit{Nonlinear Anal. }{\bf 213}:112505, 2021]. These works, when properly analyzed, leave open room for some improvement of their results. We generalize the outcomes of the mentioned articles, establish other statements and put all the claims together; in particular, we select the sharpest ones and schematize them. Moreover, we complement our research also when logistic sources are considered in the overall study. 
\end{abstract}

\section{Preamble}
For details and discussions on the meaning of the forthcoming model, especially in the frame of chemotaxis phenomena and related variants, as well as for mathematical motivations and connected state of the art, we refer to \cite{FrassuCorViglialoro,frassuviglialoro}. These articles will be often cited throughout this work.
\section{Presentation of the Theorems}\label{IntroSection}
Let $\Omega \subset \R^n$, $n \geq 2$, be a bounded and smooth domain, $\chi,\xi,\delta>0$, $m_1,m_2,m_3\in\R$, $f(u), g(u)$ and $h(u)$ be reasonably regular functions generalizing the prototypes $f(u)=K u^\alpha$, $g(u)=\gamma u^l$, and $h(u)=k u - \mu u^{\beta}$ with $K,\gamma, \mu>0$, $k \in \R$ and suitable $\alpha, l, \beta>0$. Once nonnegative initial configurations $u_0$ and $v_0$ are fixed, we aim at deriving sufficient conditions involving the above data so to ensure that the following attraction-repulsion chemotaxis model
\begin{equation}\label{problem}
\begin{cases}
u_t= \nabla \cdot ((u+1)^{m_1-1}\nabla u - \chi u(u+1)^{m_2-1}\nabla v+\xi u(u+1)^{m_3-1}\nabla w) + h(u) & \text{ in } \Omega \times (0,T_{\rm{max}}),\\
v_t=\Delta v-f(u)v  & \text{ in } \Omega \times (0,T_{\rm{max}}),\\
0= \Delta w - \delta w + g(u)& \text{ in } \Omega \times (0,T_{\rm{max}}),\\
u_{\nu}=v_{\nu}=w_{\nu}=0 & \text{ on } \partial \Omega \times (0,T_{\rm{max}}),\\
u(x,0)=u_0(x), \; v(x,0)=v_0(x) & x \in \bar\Omega,
\end{cases}
\end{equation}
admits classical solutions which are global and uniformly bounded in time. Specifically, we look for nonnegative functions $u=u(x,t), v=v(x,t), w=w(x,t)$ defined for $(x,t) \in \bar{\Omega}\times [0,T_{\rm{max}})$, and $T_{\rm{max}}=\infty$, with the properties that 
\begin{equation}\label{ClassicalAndGlobability}
\begin{cases}
u,v\in C^0(\bar{\Omega}\times [0,\infty))\cap  C^{2,1}(\bar{\Omega}\times (0,\infty)), w\in C^0(\bar{\Omega}\times [0,\infty))\cap  C^{2,0}(\bar{\Omega}\times (0,\infty)),\\ (u,v,w) \in (L^\infty((0, \infty);L^{\infty}(\Omega)))^3,
\end{cases}
\end{equation}
and pointwisely satisfying all the relations in problem \eqref{problem}.

To this scope, let  $f$, $g$ and $h$ be such that 
\begin{equation}\label{f}
f,g \in C^1(\R) \quad \textrm{with} \quad   0\leq f(s)\leq Ks^{\alpha}  \textrm{ and } \gamma s^l\leq g(s)\leq \gamma s(s+1)^{l-1},\;  \textrm{for some}\; K,\,\gamma,\,\alpha>0,\, l\geq 1 \quad \textrm{and all } s \geq 0,
\end{equation}
and 
\begin{equation}\label{h}
h \in C^1(\R) \quad \textrm{with} \quad h(0)\geq 0  \textrm{ and } h(s)\leq k s-\mu s^{\beta}, \quad  \textrm{for some}\quad k \in \R,\,\mu>0,\, \beta>1\, \quad \textrm{and all } s \geq 0.
\end{equation}
Then we prove these two theorems.
\begin{theorem}\label{MainTheorem}
Let $\Omega$ be a smooth and bounded domain of $\mathbb{R}^n$, with $n\geq 2$, $\chi, \xi, \delta$ positive reals, $l \geq 1$ and $h \equiv 0$. Moreover, for $\alpha >0$ and $m_1, m_2, m_3 \in \R$, let $f$ and $g$ fulfill \eqref{f} for each of the following cases:
\begin{enumerate}[label=$A_{\roman*}$)]
\item \label{A1} $\alpha \in \left(0, \frac{1}{n}\right]$ and $m_1>\min\left\{2m_2-(m_3+l),\max\left\{2m_2-1,\frac{n-2}{n}\right\}, m_2 - \frac{1}{n}\right\}=:\mathcal{A}$,
\item \label{A2} $\alpha \in \left(\frac{1}{n},\frac{2}{n}\right)$ and $m_1>m_2 + \alpha - \frac{2}{n}=:\mathcal{B}$,
\item \label{A3} $\alpha \in \left[\frac{2}{n},1\right]$ and $m_1>m_2 + \frac{n\alpha-2}{n\alpha-1}=:\mathcal{C}$.
\end{enumerate}
Then for any initial data $(u_0,v_0)\in (W^{1,\infty}(\Omega))^2$, with $u_0, v_0\geq 0$ on $\bar{\Omega}$, there exists a unique triplet $(u,v,w)$ of nonnegative functions, uniformly bounded in time and classically solving problem \eqref{problem}.
\end{theorem}
\begin{theorem}\label{MainTheorem1}
Under the same hypotheses of Theorem \ref{MainTheorem} and $\beta>1$, let $h$ comply with \eqref{h}.
Moreover, for $\alpha >0$ and $m_1, m_2, m_3 \in \R$, let $f$ and $g$ fulfill \eqref{f} for each of the following cases:
\begin{enumerate}[label=$A_{\roman*}$)]
\setcounter{enumi}{3}
\item \label{A4} $\alpha \in \left(0, \frac{1}{n}\right]$ and $m_1>\min\left\{2m_2-(m_3+l), 2m_2-\beta\right\}=:\mathcal{D}$,
\item \label{A5} $\alpha \in \left(\frac{1}{n},1\right)$ and $m_1>\min\left\{2m_2+1-(m_3+l), 2m_2+1-\beta\right\}=:\mathcal{E}$.
\end{enumerate}
Then the same claim holds true.
\end{theorem}
When the logistic term $h$ does not take part in the model, problem \eqref{problem} has been already analyzed in \cite{FrassuCorViglialoro} for the nonlinear diffusion and sensitivities case, and in \cite{frassuviglialoro} for the linear scenario; nevertheless, in these papers only small values of $\alpha$ are considered. Precisely, for $\alpha$ belonging to $(0,\frac{1}{2}+\frac{1}{n}),$  boundedness is ensured:
\begin{itemize}
\item in \cite[Theorem 2.1]{FrassuCorViglialoro} for $m_1,m_2,m_3\in \R$ and $l\geq 1$, under the assumption
$$m_1>\min\left\{2m_2+1-(m_3+l),\max\left\{2m_2,\frac{n-2}{n}\right\}\right\}=:\mathcal{F};$$
\item in \cite[Theorem 2.1]{frassuviglialoro} for either $m_1=m_2=m_3=l=1$, under the assumption
$$\xi>\left(\frac{8}{n}
\frac{2^\frac{2}{n}\frac{2}{n}^{n+1}(\frac{2}{n}-1)(n^2+n)}{(\frac{2}{n}+1)^{\frac{2}{n}+1}}\right)^\frac{2}{n}\lVert v_0\rVert_{L^\infty(\Omega)}^\frac{4}{n}=:\mathcal{G},$$ 
or in \cite[Theorem 2.2]{frassuviglialoro} for $m_1=m_2=m_3=1$ and any $l>1$.
\end{itemize}
In light of Theorems \ref{MainTheorem} and \ref{MainTheorem1}, herein we develop an analysis dealing also with values of $\alpha$ larger than $\frac{1}{2}+\frac{1}{n}$. Additionally, for $\alpha$ belonging to some sub-intervals of $(0,\frac{1}{2}+\frac{1}{n})$ we improve \cite[Theorem 2.1]{FrassuCorViglialoro} and  \cite[Theorems 2.1 and 2.2]{frassuviglialoro}. On the other hand, the introduction of $h$ allows us to obtain further generalizations and/or claims.

All this aspects are put together into Table \ref{Table_ResultUnified}. It, when possible, also indicates which of the assumptions taken from \cite[Theorem 2.1]{FrassuCorViglialoro}, \cite[Theorems 2.1 and 2.2]{frassuviglialoro}, and Theorems \ref{MainTheorem} and \ref{MainTheorem1} are the mildest leading to boundedness.

\setlength\extrarowheight{4.8pt}
\begin{table}[h!]
\makegapedcells
\centering 
  \begin{tabular}{ccccccccccc}
    \hline
 $m_2$ &$m_3$&$l$&\multicolumn{1}{c|}{$\alpha$}&\multicolumn{1}{c|}{$m_1$}&$\chi$&$\xi$& Reference & \quad \quad Implication\\
\hline
   $1$ &$1$&$1$&\multicolumn{1}{c|}{$[\frac{2}{n},1)$}&\multicolumn{1}{c|}{$1$}&$\R^+$&\multicolumn{1}{c|}{$ >\mathcal{G}$}& Remark \ref{Alpha}, generalizing \cite[Th. 2.1]{frassuviglialoro}\\
   $1$ &$1$&$>1$&\multicolumn{1}{c|}{$ [\frac{2}{n},1)$}&\multicolumn{1}{c|}{$1$}&$ \R^+$&\multicolumn{1}{c|}{$\R^+$}& Remark \ref{Alpha}, generalizing \cite[Th. 2.2]{frassuviglialoro}\\
    $ \R$ &$\R$&$\geq 1$&\multicolumn{1}{c|}{$ (\frac{1}{n},1)$}&\multicolumn{1}{c|}{$>\mathcal{F}$}&$\R^+$&\multicolumn{1}{c|}{$\R^+$}& Remark \ref{Alpha}, generalizing \cite[Th. 2.1]{FrassuCorViglialoro}\\
   $ \R$ &$\R$&$\geq 1$&\multicolumn{1}{c|}{$ (0,\frac{1}{n}]$}&\multicolumn{1}{c|}{$>\mathcal{A}$}&$\R^+$&\multicolumn{1}{c|}{$\R^+$}&Th. \ref{MainTheorem}& \quad \quad $\ast$ and $\ast\ast$\\
   $\R$&  $\R$&$\geq 1$&\multicolumn{1}{c|}{$(\frac{1}{n},\frac{2}{n})$}&\multicolumn{1}{c|}{$>\mathcal{B}$}&$\R^+$&\multicolumn{1}{c|}{$\R^+$}&Th. \ref{MainTheorem}& \quad \quad $\ast$\\
    $\R$&  $\R$&$\geq 1$&\multicolumn{1}{c|}{$[\frac{2}{n},1]$}&\multicolumn{1}{c|}{$>\mathcal{C}$}&$\R^+$&\multicolumn{1}{c|}{$\R^+$}&Th. \ref{MainTheorem}\\
 \hline
  $m_2$ &$m_3$&$l$&$\beta$&\multicolumn{1}{c|}{$\alpha$}&\multicolumn{1}{c|}{$m_1$}&$\chi$&$\xi$&$k$&$\mu$& Reference \\
  \hline
      $1$ &$1$&$1$&$>2$&\multicolumn{1}{c|}{$(\frac{1}{n},1)$}&\multicolumn{1}{c|}{$1$}&$\R^+$&$\R^+$&$ \R$&\multicolumn{1}{c|}{$\R^+$}&Th. \ref{MainTheorem1}\\
   $ \R$ &$\R$&$\geq 1$&$>1$&\multicolumn{1}{c|}{$ (0,\frac{1}{n}]$}&\multicolumn{1}{c|}{$>\mathcal{D}$}&$\R^+$&$ \R^+$&$ \R$&\multicolumn{1}{c|}{$\R^+$}&Th. \ref{MainTheorem1}\\
   $\R$&  $\R$&$\geq 1$&$>1$&\multicolumn{1}{c|}{$(\frac{1}{n},1)$}&\multicolumn{1}{c|}{$>\mathcal{E}$}&$\R^+$&$ \R^+$&$ \R$&\multicolumn{1}{c|}{$\R^+$}&Th. \ref{MainTheorem1}\\
   \hline  
   \end{tabular}
    \caption{Schematization collecting the ranges of the parameters involved in model \eqref{problem} for which boundedness of its solutions is established for any fixed initial distribution $u_0$ and $v_0$. The symbol $\ast$ stands for ``improves  \cite[Th. 2.1]{frassuviglialoro} and recovers \cite[Th. 2.2]{frassuviglialoro}''  and $\ast\ast$ for ``improves \cite[Th. 2.1]{FrassuCorViglialoro}''. ($\mathcal{A}, \mathcal{B}, \mathcal{C}, \mathcal{D}, \mathcal{E}, \mathcal{F}$ are defined above.)
 }\label{Table_ResultUnified}
\end{table}

\section{Local well posedness, boundedness criterion, main estimates and analysis of parameters}\label{SectionLocalInTime}
For $\Omega$, $\chi,\xi,\delta$, $m_1,m_2,m_3$ and $f, g, h$ as above, from now on with  $u, v, w \geq 0$ we refer to functions of $(x,t) \in \bar{\Omega}\times [0,T_{\rm{max}})$, for some finite $T_{\rm{max}}$, classically solving problem \eqref{problem} when nonnegative initial data $(u_0,v_0)\in (W^{1,\infty}(\Omega))^2$ are provided.  In particular, $u$ satisfies 
\begin{equation}\label{massConservation}
\int_\Omega u(x, t)dx \leq m_0 \quad \textrm{for all }\, t \in (0,\TM), 
\end{equation}
whilst $v$ is such that 
\begin{equation*}\label{MaxPrincV}
0 \leq v\leq \lVert v_0\rVert_{L^\infty(\Omega)}\quad \textrm{in}\quad \Omega \times (0,T_{\rm{max}}).
\end{equation*}
Further, globality and boundedness of $(u,v,w)$ (in  the sense of \eqref{ClassicalAndGlobability}) are ensured whenever (boundedness criterion) the $u$-component belongs to $L^\infty((0,T_{\rm{max}});L^p(\Omega))$, with $p>1$ arbitrarily large, and uniformly with respect $t\in (0,T_{\rm{max}})$. 

These basic statements can be proved by standard reasoning; in particular, when $h\equiv 0$ they verbatim follow from \cite[Lemmas 4.1 and 4.2]{FrassuCorViglialoro} and relation \eqref{massConservation} is the well-known mass conservation property. Conversely, in the presence of the logistic terms $h$ as in \eqref{h}, some straightforward adjustments have to be considered and the $L^1$-bound of $u$ is consequence of an integration of the first equation in \eqref{problem} and an application of the H\"{o}lder inequality: precisely for $k_+=\max\{k,0\}$
\[
\frac{d}{dt} \int_\Omega u = \int_\Omega h(u) =k \int_\Omega u - \mu \int_\Omega u^{\beta} \leq k_+ \int_\Omega u - \frac{\mu}{|\Omega|^{\beta-1}} \left(\int_\Omega u\right)^{\beta}
\quad \textrm{for all }\, t \in (0,\TM),
\]
and we can conclude by invoking an ODI-comparison argument.

In our computations, beyond the above positions, some uniform bounds of $\|v(\cdot,t)\|_{W^{1,s}(\Omega)}$ are required. In this sense, the following lemma gets 
the most out from $L^p$-$L^q$ (parabolic) maximal regularity; this is a cornerstone and for some small values of $\alpha$ the succeeding $W^{1,s}$-estimates are sharper than the $W^{1,2}$-estimates derived in \cite{FrassuCorViglialoro,frassuviglialoro}, and therein employed.

\begin{lemma}\label{LocalV}   
There exists $c_0>0$ such that $v$ fulfills
\begin{equation}\label{Cg}
\int_\Omega |\nabla v(\cdot, t)|^s\leq c_0 \quad \textrm{on } \,  (0,\TM)
\begin{cases}
\; \textrm{for all } s \in [1,\infty) & \textrm{if } \alpha \in \left(0, \frac{1}{n}\right],\\
\;  \textrm{for all } s \in \left[1, \frac{n}{(n\alpha-1)}\right) & \textrm{if } \alpha \in \left(\frac{1}{n},1\right].
\end{cases}  
\end{equation}
\begin{proof}
For each $\alpha \in (0,1]$, there is $\rho >\frac{1}{2}$ such that for all $s \in \left[\frac{1}{\alpha},\frac{n}{(n\alpha-1)_+}\right)$ we have 
$\frac{1}{2}<\rho <1-\frac{n}{2}\big(\alpha-\frac{1}{s}\big)$. From $1-\rho-\frac{n}{2}\big(\alpha-\frac{1}{s}\big)>0$, the claim follows invoking properties related to the Neumann heat semigroup, exactly as done in the second part of \cite[Lemma 5.1]{FrassuCorViglialoro}.
\end{proof}
\end{lemma}

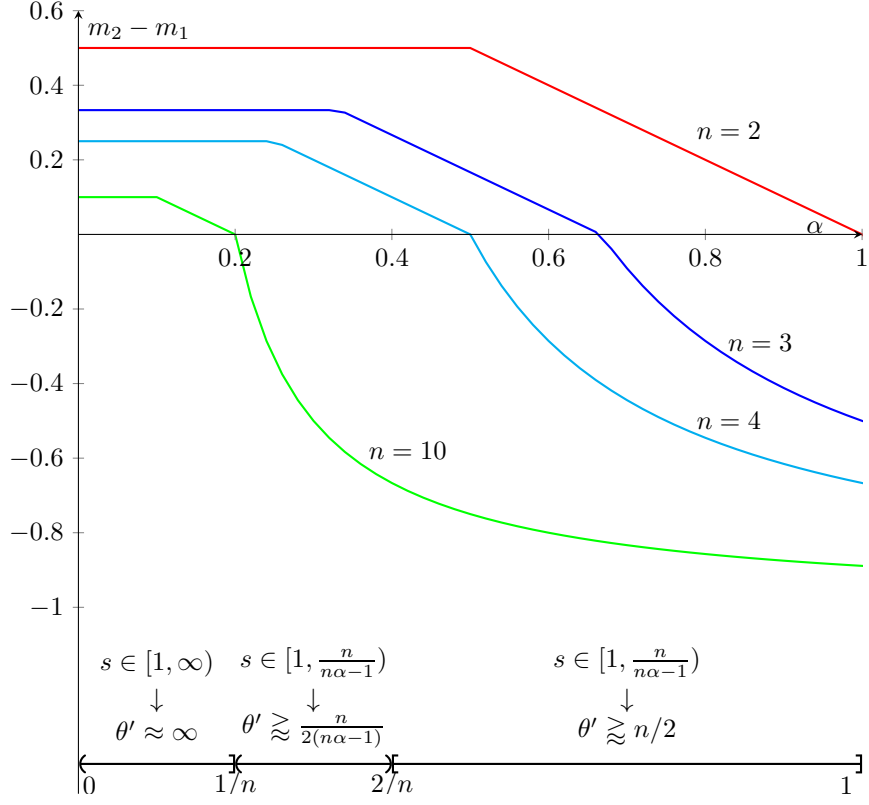
\begin{figure}
\begin{tikzpicture}[
  declare function={
    func2(\x)= (\x <= 1/2) * (1/2)   +
              and(\x > 1/2, \x <= 2/2) * (-\x+2/2)     +
              (\x >= 2/2) * (-(2*\x-2)/(2*\x-1))
   ; 
  func3(\x)= (\x <= 1/3) * (1/3)   +
              and(\x > 1/3, \x <= 2/3) * (-\x+2/3)     +
              (\x >= 2/3) * (-(3*\x-2)/(3*\x-1))
   ;
  func4(\x)= (\x <= 1/4) * (1/4)   +
              and(\x > 1/4, \x <= 2/4) * (-\x+2/4)     +
              (\x >= 2/4) * (-(4*\x-2)/(4*\x-1))
   ; 
    func5(\x)= (\x <= 1/5) * (1/5)   +
              and(\x > 1/5, \x <= 2/5) * (-\x+2/5)     +
              (\x >= 2/5) * (-(5*\x-2)/(5*\x-1))
   ; 
   func10(\x)= (\x <= 1/10) * (1/10)   +
              and(\x > 1/10, \x <= 2/10) * (-\x+2/10)     +
          (\x >= 2/10) * (-(10*\x-2)/(10*\x-1))
   ; 
  } 
]
\begin{axis}[
  axis x line=middle, axis y line=middle,
  ymin=-1.5, ymax=0.6, ytick={-1,-0.8,...,0.6}, ylabel=$m_2-m_1$,
  xmin=0, xmax=1, xtick={0,0.2,...,1},
  domain=0:2,samples=101, 
axis y line=center, axis x line=center,
			width=12cm, height=12cm,axis on top
]
\addplot [thick, red] {func2(x)};
\addplot [thick, blue] {func3(x)};
\addplot [thick, cyan] {func4(x)};
\addplot [thick, green] {func10(x)};
\draw (42,92)  node  {$n=10$};
\draw (83,100)  node  {$n=4$};
\draw (87,121)  node  {$n=3$};
\draw (83,178)  node  {$n=2$};
\draw (94,152)  node  {$\alpha$};
\draw[{(-]},  thick, black ] (0.1,8) -- (20,8);
\draw[{(-)},  thick, black] (20,8) -- (40,8);
\draw[{[-]},   thick, black ] (40,8) -- (100,8);
\draw (1.4,2.3) node   {$\tiny{0}$};
\draw (40,2.3) node   {$\tiny{2/n}$};
\draw (20,2.3) node   {$\tiny{1/n}$};
\draw (98,2.3) node   {$\tiny{1}$};
\draw (10,35) node{$\tiny{s\in [1,\infty)}$};
\draw (10,25) node{$\tiny{\downarrow}$};
\draw (10,17) node   {$\tiny{\theta' \approx \infty}$};
\draw (30,35) node  {$\tiny{s\in [1,\frac{n}{n\alpha-1})}$};
\draw (30,25) node{$\tiny{\downarrow}$};
\draw (30,17) node   {$\tiny{\theta' \gtrapprox \frac{n}{2(n\alpha-1)}}$};
\draw (70,35) node{$\tiny{s\in [1,\frac{n}{n\alpha-1})}$};
\draw (70,25) node{$\tiny{\downarrow}$};
\draw (70,17) node{$\tiny{\theta' \gtrapprox n/2}$};
\end{axis}
\end{tikzpicture}
\caption{The colored lines, functions of $\alpha$, represent the supremum of the difference $m_2-m_1$ for some space dimension $n$. Moreover, for the sub-intervals $(0,1/n]$,  $(1/n,2/n)$ and $[2/n,1]$ of $\alpha$, the corresponding range of $s$ and choice of $\theta'$ are also indicated.} \label{FigureSpiegazioneLemma}
\end{figure}

We will make use of this technical result.
\begin{lemma}\label{LemmaCoefficientAiAndExponents} 
Let $n\in \N$, with $n\geq 2$, $m_1>\frac{n-2}{n}$, $m_2,m_3\in \R$  and $\alpha \in (0,1]$. Then there is $s \in [1,\infty)$, such that for proper 
$p,q\in[1,\infty)$, $\theta$ and $\theta'$, $\mu$ and $\mu'$ conjugate exponents, we have that  
\begin{align}
a_1&= \frac{\frac{m_1+p-1}{2}\left(1-\frac{1}{(p+2m_2-m_1-1)\theta}\right)}{\frac{m_1+p-1}{2}+\frac{1}{n}-\frac{1}{2}},            
&  a_2&=\frac{q\left(\frac{1}{s}-\frac{1}{2\theta'}\right)}{\frac{q}{s}+\frac{1}{n}-\frac{1}{2}},  \nonumber \\ 
a_3 &= \frac{\frac{m_1+p-1}{2}\left(1-\frac{1}{2\alpha\mu}\right)}{\frac{m_1+p-1}{2}+\frac{1}{n}-\frac{1}{2}},  
&  a_4  &= \frac{q\left(\frac{1}{s}-\frac{1}{2(q-1)\mu'}\right)}{\frac{q}{s}+\frac{1}{n}-\frac{1}{2}},\nonumber\\ 
\kappa_1 &  =\frac{\frac{p}{2}\left(1- \frac{1}{p}\right)}{\frac{m_1+p-1}{2}+\frac{1}{n}-\frac{1}{2}} \nonumber, & \kappa_2 & =  \frac{q - \frac{1}{2}}{q+\frac{1}{n}-\frac{1}{2}}, 
\end{align}
belong to the interval $(0,1)$. If, additionally, 
\begin{equation}\label{Restrizionem1-m2-alphaPiccolo}
\alpha \in \left(0,\frac{1}{n}\right] \; \text{and}\; m_2-m_1<\frac{1}{n},
\end{equation}
\begin{equation}\label{Restrizionem1-m2-alphaGrande}
\alpha \in \left(\frac{1}{n},\frac{2}{n}\right) \; \text{and}\; m_2-m_1<\frac{2}{n}-\alpha,
\end{equation}
or
\begin{equation}\label{Restrizionem1-m2-alphaGrandeBis}
\alpha \in \left[\frac{2}{n},1\right] \; \text{and}\; m_2-m_1<\frac{2-n\alpha}{n\alpha-1},
\end{equation}
these futher relations hold true: 
\begin{equation*}\label{MainInequalityExponents}
\beta_1 + \gamma_1 =\frac{p+2m_2-m_1-1}{m_1+p-1}a_1+\frac{1}{q}a_2\in (0,1) \;\textrm{ and }\;	\beta_2 + \gamma_2= \frac{2 \alpha }{m_1+p-1}a_3+\frac{q-1}{q}a_4 
\in (0,1).
\end{equation*}
\begin{proof}
For any $s\geq 1$, let $\theta'>\max\left\{\frac{n}{2},\frac{s}{2}\right\}$  and $\mu>\max\left\{\frac{1}{2\alpha},\frac{n}{2}\right\}$. Thereafter, for 
\begin{equation}\label{Prt_q}
\begin{cases}
q > \max \left\{\frac{(n-2)}{n}\theta', \frac{s}{2\mu'}+1\right\} \\
p>\max \left\{2-\frac{2}{n}-m_1,\frac{1}{\theta}-2m_2+m_1+1, \frac{(2m_2-m_1-1)(n-2)\theta-nm_1+n}{n-(n-2)\theta},\frac{2\alpha \mu(n-2)}{n} -m_1+1\right\}, 
\end{cases}
\end{equation} 
it can be seen that $a_i,k_2\in (0,1)$, for any $i=1,2,3,4.$ On the other hand, $k_1\in (0,1)$ also thanks to the assumption $m_1>\frac{n-2}{n}.$ 

As to the second part, we distinguish three cases: $\alpha \in \left(0,\frac{1}{n}\right]$, $\alpha \in \left(\frac{1}{n},\frac{2}{n}\right)$ and 
$\alpha \in \left[\frac{2}{n},1\right].$ (We insert Figure \ref{FigureSpiegazioneLemma} to clarify the proof, by focusing on the relation involving the values of $\alpha$, $s$ and 
$\theta'$ in terms of assumptions \eqref{Restrizionem1-m2-alphaPiccolo}, \eqref{Restrizionem1-m2-alphaGrande}, \eqref{Restrizionem1-m2-alphaGrandeBis}.)
\begin{itemize}
\item [$\circ$] $\alpha \in \left(0,\frac{1}{n}\right]$. 
For $s>\frac{2\mu'}{2\mu'-1}$ arbitrarily large, consistently with \eqref{Prt_q}, we take $p=q=s$ and $\theta'=s\omega$, for some $\omega>\frac{1}{2}$. Computations provide
\begin{equation*}
0<	\beta_1+\gamma_1=\frac{s+2m_2-m_1-1-\frac{1}{\theta}}{m_1+s-2+\frac{2}{n}}+\frac{2-\frac{1}{\omega}}{s+\frac{2s}{n}},
\end{equation*}
and
\begin{equation*}
0<	\beta_2+\gamma_2=\frac{2\alpha-\frac{1}{\mu}}{m_1+s-2+\frac{2}{n}}+\frac{2s-2-\frac{s}{\mu'}}{s+\frac{2s}{n}}.
\end{equation*}
In light of the above positions, the largeness of $s$ infers $\theta$ arbitrarily close to $1$. Further, by choosing $\omega$ approaching $\frac{1}{2}$, continuity arguments imply that  $\beta_1+\gamma_1<1$ whenever restriction \eqref{Restrizionem1-m2-alphaPiccolo} is satisfied, whereas $\beta_2+\gamma_2<1$ comes from $\mu>\frac{n}{2}.$
\item  [$\circ$] $\alpha \in \left(\frac{1}{n},\frac{2}{n}\right).$ First let $s$ be arbitrarily close to $\frac{n}{n\alpha-1}$ and let $q=\frac{p}{2}$ fulfill \eqref{Prt_q}. Then, in these circumstances it holds that $\max\left\{\frac{s}{2},\frac{n}{2}\right\}=\frac{s}{2}$, so that restriction on $\theta'$ reads $\theta'>\frac{s}{2}$.  Subsequently,
\begin{equation*}
0<	\beta_1+\gamma_1=\frac{p+2m_2-m_1-1-\frac{1}{\theta}}{m_1+p-2+\frac{2}{n}}+\frac{2-\frac{s}{\theta'}}{p+\frac{2s}{n}-s},
\end{equation*}
and
\begin{equation*}
0<	\beta_2+\gamma_2=\frac{2\alpha-\frac{1}{\mu}}{m_1+p-2+\frac{2}{n}}+\frac{p-2-\frac{s}{\mu'}}{p+\frac{2s}{n}-s}.
\end{equation*}
Since from $\theta'>\frac{s}{2}$ we have that $\theta'$ approaches $\frac{n}{2(n\alpha-1)}$, similar arguments used above imply that upon enlarging $p$ one can see that condition \eqref{Restrizionem1-m2-alphaGrande} leads to  $\beta_1+\gamma_1<1$. On the other hand,  in order to have $\beta_2+\gamma_2<1$ we have to invoke the above constrain on $\mu$, i.e., $\mu>\frac{1}{2\alpha}.$
\item  [$\circ$] $\alpha \in \left[\frac{2}{n},1\right].$ We only have to consider in the previous case that now $\theta'>\frac{n}{2}$, so concluding thanks to
\eqref{Restrizionem1-m2-alphaGrandeBis}.
\end{itemize}
\end{proof}
\end{lemma}

\begin{remark}\label{RemarkOnS}
In view of its importance in the computations, we have to point out that from the above lemma $s$ can be chosen arbitrarily large only when $\alpha \in \left(0,\frac{1}{n}\right]$.
In particular, as we will see, in such an interval the terms  $\int_{\Omega} (u+1)^{p+2m_2-m_1-1} \vert \nabla v\lvert^2$ and $\int_{\Omega} (u+1)^{2\alpha} \vert \nabla v\lvert^{2(q-1)}$, appearing in our reasoning, can be treated in two alternative ways: either invoking the Young inequality or the Gagliardo--Nirenberg one. 
\end{remark}

\section{A priori estimates and proof of the Theorems}\label{EstimatesAndProofSection}

\subsection{The non-logistic case}\label{NonLog}

Recalling the globality criterion mentioned in $\S$\ref{SectionLocalInTime}, let us define the functional $y(t):=\int_\Omega (u+1)^p + \int_\Omega |\nabla v|^{2q}$, 
with $p,q>1$ properly large (and, when required, with $p=q$), and let us dedicate to derive the desired uniform bound of $\int_\Omega u^p$. 

In the spirit of Remark \ref{RemarkOnS}, let us start by analyzing the evolution in time of the functional $y(t)$ by relying on the Young inequality.
\begin{lemma}\label{Estim_general_For_u^p_nablav^2qLemma} 
Let $\alpha \in \left(0,\frac{1}{n}\right]$. If $m_1,m_2,m_3 \in \R$ comply with either $m_1>2m_2-(m_3+l)$ or $m_1>\max\left\{2m_2-1,\frac{n-2}{n}\right\}$,
then there exist $p, q>1$ such that $(u,v,w)$ satisfies for some $c_{16}, c_{17}, c_{18}>0$
\begin{equation}\label{MainInequality}
\frac{d}{dt} \left(\int_\Omega (u+1)^p + \int_\Omega  |\nabla v|^{2q}\right) + c_{16} \int_\Omega |\nabla |\nabla v|^q|^2 + c_{17} \int_\Omega |\nabla (u+1)^{\frac{m_1+p-1}{2}}|^2 \leq c_{18} \quad \text{on } (0,\TM).
\end{equation}
\begin{proof}
Let $p=q>1$ sufficiently large; moreover, in view of Remark \ref{RemarkOnS}, from now on, when necessary we will tacitly enlarge these parameters.

In the first part of the proof we focus on the estimate of the term $\frac{d}{dt} \int_\Omega (u+1)^p$.
Standard testing procedures provide 
\begin{equation*}
\begin{split}
\frac{d}{dt} \int_\Omega (u+1)^p =\int_\Omega p(u+1)^{p-1}u_t &= -p(p-1) \int_\Omega (u+1)^{p+m_1-3} |\nabla u|^2 
+p(p-1)\chi \int_\Omega u(u+1)^{m_2+p-3} \nabla u \cdot \nabla v \\
&-p(p-1)\xi \int_\Omega u(u+1)^{m_3+p-3} \nabla u \cdot \nabla w \quad \text{on } (0,\TM).
\end{split}
\end{equation*}
By reasoning as in \cite[Lemma 5.2]{FrassuCorViglialoro}, we obtain for $\epsilon_1, \epsilon_2, \tilde{\sigma}$ positive, and for all $t \in (0,\TM)$, some $c_1>0$ such that
\begin{equation}\label{Estim_1_For_u^p}  
\begin{split}
\frac{d}{dt} \int_\Omega (u+1)^p & \leq -p(p-1) \int_\Omega (u+1)^{p+m_1-3} |\nabla u|^2 +p(p-1)\chi \int_\Omega u(u+1)^{m_2+p-3} \nabla u \cdot \nabla v \\
&+ \left(\epsilon_1 + \tilde{\sigma} -  \frac{p(p-1)\xi \gamma}{2^{p-1}(m_3+p-1)} \right) \int_\Omega (u+1)^{m_3+p+l-1}
+ \epsilon_2 \int_\Omega |\nabla (u+1)^\frac{m_1+p-1}{2}|^2 + c_1.
\end{split}
\end{equation}
Let us now discuss the cases $m_1>2m_2-(m_3+l)$ and $m_1>\max\left\{2m_2-1,\frac{n-2}{n}\right\}$, respectively. 	
A double application of the Young inequality in \eqref{Estim_1_For_u^p} and bound \eqref{Cg} give
\begin{equation}\label{Young}  
\begin{split}
&p(p-1)\chi \int_\Omega u (u+1)^{m_2+p-3} \nabla u \cdot \nabla v \leq \epsilon_3 \int_\Omega (u+1)^{p+m_1-3} |\nabla u|^2 
+ c_2 \int_\Omega (u+1)^{p+2m_2-m_1-1} |\nabla v|^2\\
& \leq \epsilon_3 \int_\Omega (u+1)^{p+m_1-3} |\nabla u|^2 + \epsilon_4 \int_\Omega |\nabla v|^{s} + c_3 \int_\Omega (u+1)^{\frac{(p+2m_2-m_1-1)s}{s-2}}\\
& \leq \epsilon_3 \int_\Omega (u+1)^{p+m_1-3} |\nabla u|^2 + c_3 \int_\Omega (u+1)^{\frac{(p+2m_2-m_1-1)s}{s-2}} + c_4 \quad \text{on }\, (0,\TM),
\end{split}
\end{equation}
with $\epsilon_3, \epsilon_4 >0$ and some positive $c_2, c_3, c_4$. 
From $m_1> 2m_2-(m_3+l)$, we have $\frac{(p+2m_2-m_1-1)s}{s-2} < (m_3+p+l-1)$, and for every $\epsilon_5>0$, Young's inequality yields some $c_5>0$ entailing 
\begin{equation}\label{Young1}  
c_3 \int_\Omega (u+1)^{\frac{(p+2m_2-m_1-1)s}{s-2}} \leq \epsilon_5 \int_\Omega (u+1)^{m_3+p+l-1} + c_5 \quad \text{for all }\, (0,\TM).
\end{equation}
Now, we note that $m_1>2 m_2-1$ implies $\frac{(p+2m_2-m_1-1)s}{s-2} < p$, and the Young inequality allows us to rephrase \eqref{Young1} in an alternative way:
\begin{equation}\label{Young4}  
c_3 \int_\Omega (u+1)^{\frac{(p+2m_2-m_1-1)s}{s-2}} \leq \epsilon_6 \int_\Omega (u+1)^p + c_6 \quad \text{on }\, (0,\TM),
\end{equation}
with $\epsilon_6>0$ and positive $c_6$.
Further, an application of the Gagliardo--Nirenberg inequality and property \eqref{massConservation} yield 
\begin{equation*}\label{Theta_2}
\theta=\frac{\frac{n(m_1+p-1)}{2}\left(1-\frac{1}{p}\right)}{1-\frac{n}{2}+\frac{n(m_1+p-1)}{2}}\in (0,1),
\end{equation*}
so giving for $c_7, c_8>0$
\begin{equation*}
\begin{split}
\int_{\Omega} (u+1)^p&= \|(u+1)^{\frac{m_1+p-1}{2}}\|_{L^{\frac{2p}{m_1+p-1}}(\Omega)}^{\frac{2p}{m_1+p-1}}\\ 
&\leq c_7 \|\nabla (u+1)^{\frac{m_1+p-1}{2}}\|_{L^2(\Omega)}^{\frac{2p}{m_1+p-1}\theta} \|(u+1)^{\frac{m_1+p-1}{2}}\|_{L^{\frac{2}{m_1+p-1}}(\Omega)}^{\frac{2p}{m_1+p-1}(1-\theta)} + c_7 \|(u+1)^{\frac{m_1+p-1}{2}}\|_{L^{\frac{2}{m_1+p-1}}(\Omega)}^{\frac{2p}{m_1+p-1}}\\
& \leq c_8 \Big(\int_\Omega |\nabla (u+1)^\frac{m_1+p-1}{2}|^2\Big)^{\kappa_1}+ c_8 \quad \text{ for all } t \in(0,\TM).
\end{split}
\end{equation*}
Since $\kappa_1 \in (0,1)$ (see Lemma \ref{LemmaCoefficientAiAndExponents}), for any positive $\epsilon_7$  thanks to the Young inequality we arrive for some positive $c_9>0$ at 
 \begin{equation}\label{GN2}
\begin{split}
\epsilon_6 \int_\Omega (u+1)^p  &\leq \epsilon_7 \int_\Omega |\nabla (u+1)^\frac{m_1+p-1}{2}|^2  + c_9 \quad \text{on } (0,\TM).
\end{split}
\end{equation}
By plugging estimate \eqref{Young} into relation \eqref{Estim_1_For_u^p}, and by relying on bound \eqref{Young1} (or, alternatively to \eqref{Young1}, relations \eqref{Young4} and \eqref{GN2}), infer for appropriate $\tilde{\epsilon}_1, \tilde{\epsilon}_2>0$ and proper $c_{10}>0$ 
\begin{equation}\label{ClaimU}
\begin{split}
\frac{d}{dt} \int_\Omega (u+1)^p &\leq \left(-\frac{4p(p-1)}{(m_1+p-1)^2} + \tilde{\epsilon}_1 \right) \int_\Omega |\nabla (u+1)^{\frac{m_1+p-1}{2}}|^2\\
&+ \left(\tilde{\epsilon}_2 - \frac{p(p-1)\xi \gamma}{2^{p-1}(m_3+p-1)}\right) \int_\Omega (u+1)^{m_3+p+l-1} + c_{10} \quad \text{for all } t \in (0,\TM),
\end{split}
\end{equation}
where we also exploited that 
\begin{equation}\label{GradU}
\int_\Omega (u+1)^{p+m_1-3} |\nabla u|^2 = \frac{4}{(m_1+p-1)^2} \int_\Omega |\nabla (u+1)^{\frac{m_1+p-1}{2}}|^2 \quad  \text{on } (0,\TM).
\end{equation}
Now, as to the term $\frac{d}{dt} \int_\Omega  |\nabla v|^{2q}$ of the functional $y(t)$, reasoning similarly as in \cite[Lemma 5.3]{FrassuCorViglialoro}, we obtain for some 
$c_{11}, c_{12}>0$
\begin{equation}\label{Estim_gradV}
\frac{d}{dt}\int_\Omega  |\nabla v|^{2q}+ q \int_\Omega |\nabla v|^{2q-2} |D^2v|^2 \leq c_{11} \int_\Omega u^{2\alpha} |\nabla v|^{2q-2} +c_{12} \quad \textrm{on } \;  (0,\TM).
\end{equation}
Moreover, Young's inequalities and bound \eqref{Cg} give for every arbitrary $\epsilon_8, \epsilon_9>0$ and some $c_{13}, c_{14}, c_{15}>0$
\begin{equation}\label{Estimat_nablav^2q+2}
\begin{split}
&c_{11} \int_\Omega u^{2\alpha} |\nabla v|^{2q-2} \leq \epsilon_8 \int_\Omega (u+1)^{m_3+p+l-1} + c_{13} \int_\Omega |\nabla v|^{\frac{2(q-1)(m_3+p+l-1)}{m_3+p+l-1-2\alpha}}\\
& \leq \epsilon_8 \int_\Omega (u+1)^{m_3+p+l-1} + \epsilon_9 \int_\Omega |\nabla v|^s + c_{14} 
\leq \epsilon_8 \int_\Omega (u+1)^{m_3+p+l-1} + c_{15} \quad \textrm{for all} \quad t \in (0,\TM).
\end{split}
\end{equation}
Therefore, by inserting relation \eqref{Estimat_nablav^2q+2} into \eqref{Estim_gradV} and adding \eqref{ClaimU}, we have the claim for a proper choice of 
$\tilde{\epsilon}_1, \tilde{\epsilon}_2, \epsilon_8$ and some positive $c_{16}, c_{17}, c_{18}$, also by taking into account the relation (see \cite[page 17]{FrassuCorViglialoro})
\begin{equation}\label{GradV}
\vert \nabla \lvert \nabla v\rvert^q\rvert^2=\frac{q^2}{4}\lvert \nabla v \rvert^{2q-4}\vert \nabla \lvert \nabla v\rvert^2\rvert^2=q^2\lvert \nabla v \rvert^{2q-4}\lvert D^2v \nabla v \rvert^2\leq q^2|\nabla v|^{2q-2} |D^2v|^2.
\end{equation}
\end{proof}
\end{lemma}
Let us now turn our attention when, as mentioned before, the Gagliardo--Nirenberg inequality is employed. In this case, we can derive information not only for 
$\alpha \in \left(0,\frac{1}{n}\right]$ but also for $\alpha \in \left(\frac{1}{n},1\right]$.
\begin{lemma}\label{Met_GN}
If $m_1,m_2\in \R$ and $\alpha>0$ are taken accordingly to \eqref{Restrizionem1-m2-alphaPiccolo}, \eqref{Restrizionem1-m2-alphaGrande}, 
\eqref{Restrizionem1-m2-alphaGrandeBis}, then there exist $p, q>1$ such that $(u,v,w)$ satisfies a similar inequality as in \eqref{MainInequality}.
\begin{proof}
For $s$, $p$ and $q$ taken accordingly to Lemma \ref{LemmaCoefficientAiAndExponents} (in particular, $p=q$ for $\alpha \in \left(0,\frac{1}{n}\right]$, and $q=\frac{p}{2}$
for $\alpha \in \left(\frac{1}{n},1\right]$), let $\theta, \theta', \mu, \mu'$, $a_1,a_2, a_3, a_4$ and $\beta_1, \beta_2, \gamma_1, \gamma_2$ be therein defined.

With a view to Lemma \ref{Estim_general_For_u^p_nablav^2qLemma}, by manipulating relation \eqref{Estim_1_For_u^p} and focusing on the first inequality in \eqref{Young} and on \eqref{Estim_gradV}, proper $\epsilon_1, \tilde{\sigma}$ lead to
\begin{equation}\label{Somma}
\begin{split}
&\frac{d}{dt} \left(\int_\Omega (u+1)^p + \int_\Omega  |\nabla v|^{2q}\right) + q \int_\Omega |\nabla v|^{2q-2} |D^2v|^2 
\leq \left(-\frac{4p(p-1)}{(m_1+p-1)^2} + \tilde{\epsilon}_1 \right) \int_\Omega |\nabla (u+1)^{\frac{m_1+p-1}{2}}|^2\\
&+  c_2 \int_\Omega (u+1)^{p+2m_2-m_1-1} |\nabla v|^2 + c_{11} \int_\Omega u^{2\alpha} |\nabla v|^{2q-2} +c_{19} \quad \textrm{on } \;  (0,\TM),
\end{split}
\end{equation} 
for some $c_{19}>0$ (we also used relation \eqref{GradU}).
In this way, we can estimate the second and third integral on the right-hand side of \eqref{Somma} by applying the H\"{o}lder inequality so to have 
\begin{equation} \label{H1}
\int_{\Omega} (u+1)^{p+2m_2-m_1-1} |\nabla v|^2 \leq  \left(\int_{\Omega} (u+1)^{(p+2m_2-m_1-1)\theta}\right)^{\frac{1}{\theta}} 
\left(\int_{\Omega} |\nabla v|^{2 \theta'}\right)^{\frac{1}{\theta'}} \quad \textrm{ on } (0, \TM),
\end{equation}
and
\begin{equation} \label{H2}
\int_{\Omega} (u+1)^{2\alpha} |\nabla v|^{2q-2} \leq 
\left(\int_{\Omega} (u+1)^{2\alpha\mu}\right)^{\frac{1}{\mu}} \left(\int_{\Omega} |\nabla v|^{2(q-1)\mu'}\right)^{\frac{1}{\mu'}}\quad  \textrm{ on }  (0,\TM).
\end{equation}
By invoking the Gagliardo--Nirenberg inequality and bound \eqref{massConservation}, we obtain for some $c_{20}, c_{21}>0$
\begin{align}\label{a1}
&\left(\int_{\Omega} (u+1)^{(p+2m_2-m_1-1)\theta}\right)^{\frac{1}{\theta}}= \|(u+1)^{\frac{m_1+p-1}{2}}\|_{L^{\frac{2(p+2m_2-m_1-1)}{m_1+p-1}\theta}(\Omega)}^{\frac{2(p+2m_2-m_1-1)}{m_1+p-1}}\\ \nonumber
& \leq c_{20} \|\nabla(u+1)^{\frac{m_1+p-1}{2}}\|_{L^2(\Omega)}^{\frac{2(p+2m_2-m_1-1)}{m_1+p-1} a_1} \|(u+1)^{\frac{m_1+p-1}{2}}\|_{L^{\frac{2}{m_1+p-1}}(\Omega)}^{\frac{2(p+2m_2-m_1-1)}{m_1+p-1} (1-a_1)} + c_{20} \|(u+1)^{\frac{m_1+p-1}{2}}\|_{L^{\frac{2}{m_1+p-1}}(\Omega)}^{\frac{2(p+2m_2-m_1-1)}{m_1+p-1}} \\ \nonumber
&\leq c_{21} \left(\int_{\Omega} |\nabla (u+1)^{\frac{m_1+p-1}{2}}|^2\right)^{\beta_1}+ c_{21} \quad \textrm{ for all } \,t\in(0,\TM),
\end{align}
and for some $c_{22}, c_{23}>0$
\begin{align}\label{a3}
&\left(\int_{\Omega} (u+1)^{2\alpha\mu}\right)^{\frac{1}{\mu}}= \|(u+1)^{\frac{m_1+p-1}{2}}\|_{L^{\frac{4\alpha\mu}{m_1+p-1}}(\Omega)}^{\frac{4\alpha}{m_1+p-1}}\\ \nonumber
& \leq c_{22} \|\nabla(u+1)^{\frac{m_1+p-1}{2}}\|_{L^2(\Omega)}^{\frac{4\alpha}{m_1+p-1}  a_3} \|(u+1)^{\frac{m_1+p-1}{2}}\|_{L^{\frac{2}{m_1+p-1}}(\Omega)}^{\frac{4\alpha}{m_1+p-1} (1-a_3)} + c_{22} \|(u+1)^{\frac{m_1+p-1}{2}}\|_{L^{\frac{2}{m_1+p-1}}(\Omega)}^{\frac{4\alpha}{m_1+p-1}}\\ \nonumber
&\leq c_{23} \left(\int_{\Omega} |\nabla (u+1)^{\frac{m_1+p-1}{2}}|^2\right)^{\beta_2}+ c_{23}
\quad \textrm{for all }\, t \in (0,\TM).
\end{align}

In a similar way, we can again apply the Gagliardo--Nirenberg inequality and bound \eqref{Cg} and get for some $c_{24}, c_{25}>0$
\begin{align}\label{a2}
&\left(\int_{\Omega} |\nabla v|^{2 \theta'}\right)^{\frac{1}{\theta'}} =\| |\nabla v|^q \|_{L^{\frac{2 \theta'}{q}}(\Omega)}^{\frac{2}{q}}\leq c_{24} \|\nabla |\nabla v|^q\|_{L^2(\Omega)}^{\frac{2}{q}a_2} \| |\nabla v|^q\|_{L^{\frac{s}{q}}(\Omega)}^{\frac{2}{q}(1-a_2)} + c_{24} \| |\nabla v|^q\|_{L^{\frac{s}{q}}(\Omega)}^{\frac{2}{q}} \\ \nonumber
& \leq c_{25} \left(\int_{\Omega} |\nabla |\nabla v|^q|^2 \right)^{\gamma_1}+ c_{25}
\quad \textrm{ for all } t\in (0, \TM),
\end{align}
and for some $c_{26}, c_{27}>0$
\begin{align}\label{a4}
&\left(\int_{\Omega} |\nabla v|^{2(q-1) \mu'}\right)^{\frac{1}{\mu'}} =\| |\nabla v|^q \|_{L^{\frac{2(q-1)}{q}\mu'}(\Omega)}^{\frac{2(q-1)}{q}} 
\leq c_{26} \|\nabla |\nabla v|^q\|_{L^2(\Omega)}^{\frac{2(q-1)}{q}  a_4} \| |\nabla v|^q\|_{L^{\frac{s}{q}}(\Omega)}^{\frac{2(q-1)}{q}  (1-a_4)}
+ c_{26} \| |\nabla v|^q\|_{L^{\frac{s}{q}}(\Omega)}^{\frac{2(q-1)}{q}} \\ \nonumber
&\leq c_{27} \left(\int_{\Omega} |\nabla |\nabla v|^q|^2 \right)^{\gamma_2}+ c_{27} \quad \textrm{ for every } t\in (0,\TM).
\end{align}

By plugging \eqref{H1} and \eqref{H2} into \eqref{Somma} and taking into account \eqref{a1}, \eqref{a3}, \eqref{a2}, \eqref{a4}, we deduce for a proper 
$\tilde{\epsilon}_1$, once inequality \eqref{GradV} is considered, the following estimate for some $c_{28}, c_{29}, c_{30}, c_{31}, c_{32}>0$:
\begin{equation}\label{Somma1}
\begin{split}
&\frac{d}{dt} \left(\int_\Omega (u+1)^p + \int_\Omega  |\nabla v|^{2q}\right) + c_{28} \int_\Omega |\nabla |\nabla v|^q|^2 
+ c_{29} \int_\Omega |\nabla (u+1)^{\frac{m_1+p-1}{2}}|^2 - c_{30}\\
&\leq c_{31} \left(\int_{\Omega} |\nabla (u+1)^{\frac{m_1+p-1}{2}}|^2\right)^{\beta_1}  
\left(\int_{\Omega} |\nabla |\nabla v|^q|^2 \right)^{\gamma_1} + c_{31} \left(\int_{\Omega} |\nabla (u+1)^{\frac{m_1+p-1}{2}}|^2\right)^{\beta_1}\\
& + c_{31} \left(\int_{\Omega} |\nabla |\nabla v|^q|^2 \right)^{\gamma_1} 
+ c_{32} \left(\int_{\Omega} |\nabla (u+1)^{\frac{m_1+p-1}{2}}|^2\right)^{\beta_2} \left(\int_{\Omega} |\nabla |\nabla v|^q|^2 \right)^{\gamma_2}\\ 
&+ c_{32} \left(\int_{\Omega} |\nabla (u+1)^{\frac{m_1+p-1}{2}}|^2\right)^{\beta_2} + c_{32} \left(\int_{\Omega} |\nabla |\nabla v|^q|^2 \right)^{\gamma_2} \quad \text{on }
(0,\TM).
\end{split}
\end{equation} 
Since by Lemma \ref{LemmaCoefficientAiAndExponents} we have that $\beta_1 + \gamma_1 <1$ and $\beta_2 + \gamma_2 <1$, and in particular $\beta_1, \gamma_1, \beta_2, \gamma_2 \in (0,1)$, we can treat the two integral products and the remaining four addenda of the right-hand side in a such way that eventually they are absorbed by the two integral terms involving the gradients in the left one. More exactly, to the products we apply 
\[
a^{d_1}b^{d_2} \leq \epsilon(a+b)+c \quad \textrm{with } a,b\geq0, d_1,d_2 >0 \; \textrm{such that } d_1+d_2<1, \; \textrm{for all } \epsilon>0 \; \textrm{and some } c>0 
\]
(achievable by means of applications of Young's inequality), and to the other terms the Young inequality. In this way, the resulting linear combination of $\int_{\Omega} |\nabla |\nabla v|^q|^2$ and $\int_{\Omega} |\nabla (u+1)^{\frac{m_1+p-1}{2}}|^2$ can be turned into $\frac{c_{28}}{2} \int_{\Omega} |\nabla |\nabla v|^q|^2 + \frac{c_{29}}{2} \int_{\Omega} |\nabla (u+1)^{\frac{m_1+p-1}{2}}|^2$, which coming back to \eqref{Somma1} infers what claimed.
\end{proof}
\end{lemma}

\subsection{The logistic case}\label{Log}
For the logistic case we retrace part of the computations above connected to the usage of the Young inequality only.
\begin{lemma}\label{Estim_general_For_u^p_nablav^2qLemmaLog}
If $m_1,m_2,m_3 \in \R$ comply with $m_1>2m_2-(m_3+l)$ or $m_1>2m_2-\beta$ whenever $\alpha \in \left(0, \frac{1}{n}\right]$, or 
$m_1>2m_2 +1-(m_3+l)$ or $m_1>2m_2+1-\beta$ whenever $\alpha \in \left(\frac{1}{n},1\right)$, then there exist $p,q>1$ such that $(u,v,w)$ satisfies 
a similar inequality as in \eqref{MainInequality}.
\begin{proof}
As in Lemma \ref{Estim_general_For_u^p_nablav^2qLemma}, in view of inequality \eqref{Young} and the properties of the logistic $h$ in \eqref{h}, relation \eqref{Estim_1_For_u^p} now becomes for some positive $c_{33}$
\begin{equation}\label{Estim_1_For_u^p1}  
\begin{split}
\frac{d}{dt} \int_\Omega (u+1)^p & \leq (-p(p-1)+\epsilon_3) \int_\Omega (u+1)^{p+m_1-3} |\nabla u|^2 +c_3 \int_\Omega (u+1)^{\frac{(p+2m_2-m_1-1)s}{s-2}}\\
&+ \left(\epsilon_1 + \tilde{\sigma} -  \frac{p(p-1)\xi \gamma}{2^{p-1}(m_3+p-1)} \right) \int_\Omega (u+1)^{m_3+p+l-1}
+ \epsilon_2 \int_\Omega |\nabla (u+1)^\frac{m_1+p-1}{2}|^2\\
&+ pk_+ \int_\Omega (u+1)^p - p \mu  \int_\Omega (u+1)^{p-1}u^{\beta} + c_{33} \quad \text{ for all } t \in (0,\TM).
\end{split}
\end{equation}
Applying the inequality $(A+B)^p \leq 2^{p-1} (A^p+B^p)$ with $A,B \geq 0$ and $p>1$ to the last integral in \eqref{Estim_1_For_u^p1}, implies that  
$-u^{\beta} \leq -\frac{1}{2^{\beta-1}} (u+1)^{\beta}+1$; therefore 
\begin{equation}\label{beta}
- p \mu  \int_\Omega (u+1)^{p-1} u^{\beta} \leq -\frac{p \mu}{2^{\beta-1}} \int_\Omega (u+1)^{p-1+\beta} + p \mu \int_\Omega (u+1)^{p-1} \quad \text{on } (0,\TM).
\end{equation}
Henceforth, by taking into account the Young inequality, we have that for $t \in (0,\TM)$
\begin{equation} \label{k}
pk_+ \int_\Omega (u+1)^p \leq \delta_1 \int_\Omega (u+1)^{p-1+\beta} + c_{34} \quad \text{and} \quad  
p \mu \int_\Omega (u+1)^{p-1} \leq \delta_2 \int_\Omega (u+1)^{p-1+\beta} + c_{35},
\end{equation}
with $\delta_1, \delta_2 >0$ and some $c_{34}, c_{35} >0$.

\quad \textbf{Case $1$}: $\alpha \in \left(0, \frac{1}{n}\right]$ and $m_1>2m_2-(m_3+l)$ or $m_1>2m_2-\beta$.
For $m_1> 2m_2-(m_3+l)$ we refer to Lemma \ref{Estim_general_For_u^p_nablav^2qLemma} and we take in mind inequality \eqref{Young1}. 
Conversely, when $m_1>2 m_2-\beta$, we have that (recall $s$ may be arbitrary large) $\frac{(p+2m_2-m_1-1)s}{s-2} < p-1+\beta$, and by means of the Young inequality estimate \eqref{Young1} can alternatively read
\begin{equation}\label{young4}  
c_3 \int_\Omega (u+1)^{\frac{(p+2m_2-m_1-1)s}{s-2}} \leq \delta_3 \int_\Omega (u+1)^{p-1+\beta} + c_{36} \quad \text{on }\, (0,\TM),
\end{equation}
with $\delta_3>0$ and positive $c_{36}$.
By inserting estimates \eqref{beta} and \eqref{k} into relation \eqref{Estim_1_For_u^p1}, as well as taking into account \eqref{Young1} 
(or, alternatively to \eqref{Young1}, bound \eqref{young4}), for suitable $\hat{\epsilon}, \tilde{\epsilon}_2, \tilde{\delta} >0$ and some $c_{37}>0$ we arrive at 
\begin{equation*}\label{ClaimUlog}
\begin{split}
\frac{d}{dt} \int_\Omega (u+1)^p &\leq \left(-\frac{4p(p-1)}{(m_1+p-1)^2} + \hat{\epsilon}\right) \int_\Omega |\nabla (u+1)^{\frac{m_1+p-1}{2}}|^2 
+ \left(\tilde{\epsilon_2} - \frac{p(p-1)\xi \gamma}{2^{p-1}(m_3+p-1)} \right) \int_\Omega (u+1)^{m_3+p+l-1}\\
& + \left(\tilde{\delta} - \frac{p \mu}{2^{\beta-1}} \right) \int_\Omega (u+1)^{p-1+\beta} + c_{37} \quad \text{ for all } t \in (0,\TM),
\end{split}
\end{equation*}
where we used again relation \eqref{GradU}. We conclude reasoning exactly as in the second part of the proof of Lemma \ref{Estim_general_For_u^p_nablav^2qLemma} and by choosing suitable 
$\hat{\epsilon}, \tilde{\epsilon}_2, \tilde{\delta}, \epsilon_8$.

\textbf{Case $2$}: $\alpha \in \left(\frac{1}{n},1\right)$ and $m_1>2m_2 +1-(m_3+l)$ or $m_1>2m_2+1-\beta$.
Accordingly to Remark \ref{RemarkOnS}, since now $s$ cannot increase arbitrarily, relations \eqref{Young1} and \eqref{young4} have to be differently discussed. 
In particular, for some $\bar{c}_1>0$ we can estimate relation \eqref{Young} as follows:
\begin{equation*}
\begin{split}
&p(p-1)\chi \int_\Omega u (u+1)^{m_2+p-3} \nabla u \cdot \nabla v \leq \epsilon_3 \int_\Omega (u+1)^{p+m_1-3} |\nabla u|^2 
+ c_2 \int_\Omega (u+1)^{p+2m_2-m_1-1} |\nabla v|^2\\
& \leq \epsilon_3 \int_\Omega (u+1)^{p+m_1-3} |\nabla u|^2 + \epsilon_4 \int_\Omega |\nabla v|^{2(p+1)} + \bar{c}_1 \int_\Omega (u+1)^{\frac{(p+2m_2-m_1-1)(p+1)}{p}} 
\quad \text{on }\, (0,\TM).
\end{split}
\end{equation*}
Now, if $m_1>2m_2+1-(m_3+l)$, then some $p$ sufficiently large infers to $\frac{(p+2m_2-m_1-1)(p+1)}{p}< p+m_3+l-1$, so that for any positive $\bar{\epsilon}_1$ 
and some $\bar{c}_2>0$ we have
\[
\bar{c}_1 \int_\Omega (u+1)^{\frac{(p+2m_2-m_1-1)(p+1)}{p}} \leq \bar{\epsilon}_1 \int_\Omega (u+1)^{p+m_3+l-1} + \bar{c}_2 \quad \text{on } (0,\TM).
\]
Conversely, and in a similar way, for $m_1>2m_2+1-\beta$ we have for any positive $\bar{\epsilon}_2$ and some $\bar{c}_3>0$
\[
\bar{c}_1 \int_\Omega (u+1)^{\frac{(p+2m_2-m_1-1)(p+1)}{p}} \leq \bar{\epsilon}_2 \int_\Omega (u+1)^{p-1+\beta} + \bar{c}_3 \quad \text{on } (0,\TM).
\]
The remaining part of the proof follows as the previous case, by taking into account \cite[Lemma 5.3-Lemma 5.4]{FrassuCorViglialoro} for the term dealing with
$ \int_\Omega |\nabla v|^{2(p+1)}$. 
\end{proof}
\end{lemma}
As a by-product of what now obtained we are in a position to conclude.
\subsection{Proof of Theorems \ref{MainTheorem} and \ref{MainTheorem1}}
\begin{proof}
Let $(u_0,v_0) \in (W^{1,\infty}(\Omega))^2$ with $u_0, v_0 \geq 0$ on $\bar{\Omega}$. For $f$ and $g$ as in \eqref{f} and, respectively, for $f$, $g$ as in \eqref{f} and
$h$ as in \eqref{h}, let $\alpha >0$ and let $m_1,m_2,m_3 \in \R$ comply with \ref{A1}, \ref{A2} and \ref{A3}, respectively, \ref{A4} and \ref{A5}. Then, we refer to Lemmas \ref{Estim_general_For_u^p_nablav^2qLemma} and \ref{Met_GN}, respectively, Lemma \ref{Estim_general_For_u^p_nablav^2qLemmaLog} and obtain for some 
$C_1,C_2,C_3>0$
\begin{equation}\label{Estim_general_For_y_2}
y'(t) + C_1 \int_\Omega |\nabla (u+1)^{\frac{m_1+p-1}{2}}|^2 + C_2 \int_\Omega |\nabla |\nabla v|^q|^2 \leq C_3
\quad \text{ on } (0, \TM).
\end{equation}
Successively, the Gagliardo--Nirenberg inequality again makes that some positive constants $c_{38}, c_{39}$ imply from the one hand
\begin{equation*}\label{G_N2}
\int_\Omega (u+1)^p \leq c_{38} \Big(\int_\Omega |\nabla (u+1)^\frac{m_1+p-1}{2}|^2\Big)^{\kappa_1}+ c_{38}  \quad \text{for all } t \in (0,\TM),
\end{equation*}
(as already done in \eqref{GN2}), and from the other
\begin{equation}\label{Estim_Nabla nabla v^q}
\int_\Omega \lvert \nabla v\rvert^{2q}=\lvert \lvert \lvert \nabla v\rvert^q\lvert \lvert_{L^2(\Omega)}^2 
\leq c_{39} \lvert \lvert\nabla  \lvert \nabla v \rvert^q\rvert \lvert_{L^2(\Omega)}^{2\kappa_2} \lvert \lvert\lvert \nabla v \rvert^q\lvert \lvert_{L^\frac{1}{q}(\Omega)}^{2(1-\kappa_2)} +c_{39} \lvert \lvert \lvert \nabla v \rvert^q\lvert \lvert^2_{L^\frac{1}{q}(\Omega)}\quad \textrm{on } (0,\TM),
\end{equation}
with $\kappa_2$ already defined in Lemma \ref{LemmaCoefficientAiAndExponents}. 
Subsequently, the $L^s$-bound of $\nabla v$ in \eqref{Cg} infers some $c_{40}>0$ such that
\begin{equation*}\label{Estim_Nabla nabla v^p^2}
\int_\Omega \lvert \nabla v\rvert^{2q}\leq c_{40} \Big(\int_\Omega \lvert \nabla \lvert \nabla v \rvert^q\rvert^2\Big)^{\kappa_2}+c_{40} \quad \text{for all } t \in (0,\TM).
\end{equation*}
In the same flavour of \cite[Lemma 5.4]{FrassuCorViglialoro}, by using the estimates involving $\int_\Omega (u+1)^p$ and $\int_\Omega \rvert\nabla v\lvert^{2q}$,
relation \eqref{Estim_general_For_y_2} provides positive constants  $c_{41}$ and $c_{42}$, and $\kappa=\min\{\frac{1}{\kappa_1},\frac{1}{\kappa_2}\}$ such that 
\begin{equation*}\label{MainInitialProblemWithM}
\begin{cases}
y'(t)\leq c_{41}-c_{42} y^{\kappa}(t)\quad \textrm{for all } t \in (0,T_{\rm{max}}),\\
y(0)=\int_\Omega (u_0+1)^p+ \int_\Omega |\nabla v_0|^{2q}. 
\end{cases}
\end{equation*}
Finally, ODE comparison principles imply $u \in L^\infty((0,\TM);L^p(\Omega))$, and the conclusion is a consequence of the boundedness criterion in 
$\S$\ref{SectionLocalInTime}.
\end{proof}

\begin{remark}[On the validity of the theorems in \cite{FrassuCorViglialoro} and \cite{frassuviglialoro} for $\alpha \geq \frac{1}{2}+\frac{1}{n}$]\label{Alpha}
In the proofs of \cite[Theorem 2.1]{FrassuCorViglialoro} and \cite[Theorems 2.1 and 2.2]{frassuviglialoro}, it is seen that the $L^2$ uniform estimate of $\nabla v$ is used to control some integral on $\partial \Omega$ (and this allows us to avoid to restrict our study to convex domains), as well as to deal with the term 
$\int_\Omega \rvert \nabla v \lvert^{2p}$ with the Gagliardo--Nirenberg inequality; for instance we are referring to \cite[(28) and (39)]{frassuviglialoro}, respectively. 
Such finiteness of $\int_\Omega \rvert \nabla v \lvert^2$ is related to the values of $\alpha$ in these articles:
$\alpha \in \left(0, \frac{1}{2}+\frac{1}{n}\right)$ (see \cite[Lemma 4.1]{frassuviglialoro}). Apparently only $\nabla v \in L^\infty((0,\TM);L^1(\Omega))$ suffices to address these issues. Indeed, as far as the topological property of $\Omega$ is concerned, we can invoke \cite[(3.10) of Proposition 8]{YokotaEtAlNonCONVEX} with $s=1$; on the other hand, for the question tied to the employment of the Gagliardo--Nirenberg inequality, we may operate as done in \eqref{Estim_Nabla nabla v^q}. As a consequence, in view of Lemma \ref{LocalV}, we have that  $\nabla v \in L^\infty((0,\TM);L^1(\Omega))$, so that \cite[Theorem 2.1]{FrassuCorViglialoro} and \cite[Theorems 2.1 and 2.2]{frassuviglialoro} hold true for any $\alpha \in (0,1)$.
\end{remark}

\subsubsection*{Acknowledgments}
SF and GV are members of the Gruppo Nazionale per l'Analisi Matematica, la Probabilit\`a e le loro Applicazioni (GNAMPA) of the Istituto Nazionale di Alta Matematica (INdAM) and are partially supported by the research project \emph{Evolutive and Stationary Partial Differential Equations with a Focus on Biomathematics}, funded by Fondazione di Sardegna (2019). GV is partially supported by MIUR (Italian Ministry of Education, University and Research) Prin 2017 \emph{Nonlinear Differential Problems via Variational, Topological and Set-valued Methods} (Grant Number: 2017AYM8XW). TL is partially supported by NNSF of P. R. China (Grant No. 61503171), CPSF (Grant No. 2015M582091), and NSF of Shandong Province (Grant No. ZR2016JL021).


\begin{thebibliography}{1}

\bibitem{FrassuCorViglialoro}
S.~Frassu, C.~{van der Mee}, and G.~Viglialoro.
\newblock Boundedness in a nonlinear attraction-repulsion {K}eller--{S}egel
  system with production and consumption.
\newblock {\em J. Math. Anal. Appl.}, 504(2):125428, 2021.

\bibitem{frassuviglialoro}
S.~Frassu and G.~Viglialoro.
\newblock Boundedness in a chemotaxis system with consumed chemoattractant and
  produced chemorepellent.
\newblock {\em Nonlinear Anal.}, 213:112505, 2021.

\bibitem{YokotaEtAlNonCONVEX}
S.~Ishida, K.~Seki, and T.~Yokota.
\newblock Boundedness in quasilinear {K}eller--{S}egel systems of
  parabolic-parabolic type on non-convex bounded domains.
\newblock {\em J. Differential Equations}, 256(8):2993--3010, 2014.

\end{thebibliography}

\end{document}